\documentclass[]{article}
\usepackage{amssymb,amsmath,amscd,amsthm,xcolor}
\usepackage{authblk}
\usepackage[cm]{fullpage}
\usepackage{graphicx}
\usepackage{multirow}
\newcommand{\Keywords}[1]{\par\noindent
{\small{\em Keywords\/}: #1}}

\definecolor{pusu}{rgb}{0.0, 0.42, 0.24}

\newcommand{\Rd}{\mathbb{R}}

\newcommand{\Un}{\mathbf{1}}

\newcommand{\eps}{\varepsilon}
\newcommand{\fy}{\varphi}

\newcommand{\tend}[1]{\longrightarrow}
\newcommand{\tends}[1]{\xrightarrow[#1]{}}
\newcommand{\tendsb}{\xrightarrow{a.s.}}
\newcommand{\tendsd}{\xrightarrow{\ d\ }}        
\newcommand{\imply}{\Rightarrow}             %

\newcommand{\norm}[1]{\lVert #1\rVert}                    

\newcommand{\abs}[1]{\lvert #1\rvert}                     
\newcommand{\Abs}[1]{\bigl| #1\bigr|}

\DeclareMathOperator{\Prob}{P}

\renewcommand{\d}{\mathrm{d}}

\let\0=\varnothing
\let\1=\Un

\renewcommand{\le}{\leqslant}
\renewcommand{\ge}{\geqslant}

\newcommand{\edoc}{\end{document}}

\let\s=\sigma

\let\l=\lambda
\let\a=\alpha

\let\de=\delta
\let\D=\Delta
\let\t=\theta
\let\ls=\leqslant
\let\gs=\geqslant
\def\bt{\hbox{$\bf \cdot$}}
\def\E{\hbox{\bf E}}
\def\pr{\hbox{\bf P}}
\let\kap=\varkappa
\def\dr{\hbox{\rm d}}
\def\bt{\hbox{$\mathbf \cdot$}}

\theoremstyle{plain}
\newtheorem{theorem}{\sc Theorem}[section]
\newtheorem{lemma}{\sc Lemma}[section]
\newtheorem{proposition}{\sc Proposition}[section]

\theoremstyle{definition}

\newtheorem{remark}{\sc Remark}[section]

\numberwithin{equation}{section}

\begin{document}
\date{}

\title{On some estimators of the Hurst index of the solution of SDE driven by a fractional Brownian motion}
\author[1,*,\dag]{K. Kubilius}
\author[2,**]{V. Skorniakov}

\affil[1]{Vilnius University, Institute of Mathematics and
Informatics, Akademijos 4, LT-08663, Vilnius, Lithuania}
\affil[2]{Vilnius University, Faculty of Mathematics and
Informatics, Naugarduko 24, LT-03225, Vilnius, Lithuania}

\maketitle

\let\oldthefootnote\thefootnote
\renewcommand{\thefootnote}{\fnsymbol{footnote}}
\footnotetext[1]{Corresponding author. E-mail:
\texttt{kestutis.kubilius@mii.vu.lt}} \footnotetext[2]{This research was funded by a grant (No. MIP-048/2014) from the Research Council of Lithuania.} \footnotetext[7]{E-mail:
\texttt{viktor.skorniakov@mif.vu.lt}}
\let\thefootnote\oldthefootnote

\abstract{Strongly consistent and asymptotically normal estimators
of the Hurst parameter of solutions of stochastic differential equations are
proposed. The estimators are based on discrete observations of the
underlying processes.
\bigskip
\Keywords{fractional
Brownian motion, stochastic differential equation, Hurst index} }

\section{Introduction}\label{s:main_result}

Recently long range dependence (LRD) became one of the most researched phenomena in statistics. It appears in various applied fields and inspires new models to account for it. Stochastic differential equations (SDEs) are widely used to model continuous time processes. Within this framework, LRD is frequently modeled with the help of SDEs driven by a fractional Brownian motion (fBm). It is well known that the latter Gaussian process is governed by a single parameter $H\in(0,1)$ (called the Hurst index) and that values of $H$ in $(1/2,1)$ correspond to LRD models. In applications,  the estimation of $H$ is a fundamental problem. Its solution depends on the theoretical structure of a model under consideration. Therefore, particular models usually deserve separate analysis. In this paper, we concentrate on the estimation of $H$ under the assumption that an observable continuous time process $(X_t)_{t\in[0,T]}$ satisfies SDE
\begin{equation}\label{e:diflygt}
X_t=\xi + \int_0^t f(X_s)\,\d s + \int_0^t g(X_s)\,\d B^H_s,\quad t\in[0,T],
\end{equation}
where $T>0$ is fixed, $\xi$ is an initial r.v., $f$ and $g$ are continuous functions satisfying some regularity conditions, and $(B^H_t)_{t\in[0,T]}$ is a fBm with the Hurst index $1/2<H<1$. Our goal is to construct a strongly consistent and asymptotically normal estimator of the $H$ from discrete observations $X_{t_1},\dots, X_{t_n}$ of trajectory $X_t$, $t\in[0,T]$.

We consider two cases. First, we examine the case when $g$ is completely specified. Next, we relax this restriction and allow $g$ to be unknown. Such situation may appear when $g$ depends on additional nuisance parameters. In both cases, boundedness of $1/g$ plays an important role and is assumed to hold.

To our best knowledge, so far only several studies investigated this question. The  pioneering work has been done by \cite{Berzin-08} of Berzin and Le\'on as well as lecture notes \cite{BLB} with references therein. \cite{kubmel1}, \cite{kubmel2} and \cite{kubmish} were also devoted to the problems of the same nature. However, all of these works focused on the strong consistency. The present paper is a generalization of \cite{ksm} where a special case of (\ref{e:diflygt}) was considered.

The paper is organized in the following way. In Section \ref{s:main_result} we present the main results of the paper. Section \ref{s:preliminaries} is devoted to several results needed for the proofs. Sections \ref{s:proofs1}--\ref{s:main} contain the proofs of the main results. Finally, in Section \ref{s:examples} two examples are given in order to illustrate the obtained results.

\section{Main results}\label{s:main_result}

To avoid cumbersome expressions, we introduce symbols
$O_{\omega},o_{\omega}$. Let $(Y_n)$ be a sequence of r.v.,
$\varsigma$ is an a.s. non-negative r.v. and
$(a_n)\subset(0,\infty)$ vanishes. $Y_n=O_\omega(a_n)$ means that
$\vert Y_n\vert\le \varsigma\cdot a_n$; $Y_n=o_\omega(a_n)$ means
that $\vert Y_n\vert\le \varsigma\cdot b_n$ with $b_n=o(a_n)$. In
particular, $Y_n=o_\omega(1)$ corresponds to the sequence $(Y_n)$
which tends to $0$ a.s. as $n\to \infty$.

Let $\pi_n=\{\tau^n_k,\ k=0,\dots,i_n\}$, $n\ge1$, $\mathbb{N}\ni
i_n\uparrow\infty$, be a sequence of partitions of the interval $[0,T]$. If  partition $\pi_n$ is uniform  then $\tau^n_k=\frac{kT}{i_n}$ for  all $k\in\{0,\ldots,i_n\}$. If $i_n\equiv n$, we write $t_k^n$ instead of $\tau_k^n$. In order
to formulate our main results, we state two hypotheses:
\begin{align}
&\mathbf{(H)}\quad  \Delta X_{\tau^n_k} =X_{\tau^n_k}-X_{\tau^n_{k-1}}=O_\omega\big(d_n^{H-\eps}\big), \quad k=1,\dots, i_n;\nonumber\\
&\mathbf{(H_1)}\quad  \Delta^{(2)}
X_{\tau^n_k}=X_{\tau^n_k}-2X_{\tau^n_{k-1}} +X_{\tau^n_{k-2}}
=g(X_{\tau^n_{k-1}})\Delta^{(2)} B^{H}_{\tau^n_k}
+O_\omega\big(d_n^{2(H-\eps)}\big),\label{e:nario_asimp_H1}\\
&\qquad \qquad k=2,\dots,i_n, \nonumber
\end{align}
for all $\eps\in(0,H-\frac 12)$, where $d_n=\max_{1\ls k\ls i_n}(\tau^n_k-\tau^n_{k-1})$.

\begin{theorem}\label{t:case} Assume that solution of Eq. (\ref{e:diflygt}) satisfies hypothesis $\mathbf{(H_1)}$. Moreover, let
the function $g$ is known, Lipschitz-continuous and there exists a
random variable $\varsigma$ such that $\Prob(\varsigma<\infty)=1$
and
\begin{equation}\label{i:boundness}
\sup_{t\in[0,T]}\frac{1}{\vert g(X_t)\vert}\le \varsigma\quad\mbox{a.s.}
\end{equation}
Then
\begin{align*}
\widehat H^{(1)}_n&\tend{} H \quad\mbox{a.s.},\\
2\sqrt{n}\,\ln\frac{n}{T}\,(\widehat H^{(1)}_n-H)&\tendsd
N(0;\sigma^2_H)\qquad\mbox{for}\ H\in(1/2,1),
\end{align*}
where  $\sigma^2_H$ is a known variance  defined in Subsection \ref{s:fbm},
\begin{align*}
\widehat H^{(1)}_n=&\fy_{n,T}^{-1}\Bigg(\frac{1}{n}
\sum_{i=2}^n\left(\frac{\Delta^{(2)}
X_{t^n_i}}{g\big(X(t^n_{i-1})\big)}\right)^2\Bigg) \quad \mbox{for}\
n>T,\\
\fy_{n,T}(x)=&\Big(\frac{T}{n}\Big)^{2x}(4-2^{2x}) \text{  and  }
\fy_{n,T}^{-1}\ \mbox{denotes an inverse of}\ \fy_{n,T},\
x\in(0,1),\ n>T.\qquad\qedsymbol
\end{align*}
\end{theorem}

It is natural to try to drop restriction of the known $g$. For this
purpose we need several additional  definitions. Assume that the
process $X$ is observed at time points $\frac{i}{m_n}T$,
$i=1,\dots,m_n$, where $m_n=n k_n$,  and $k_n$ grows faster than
$n\ln n$, but the growth does not exceed polynomial, e.g.
$k_n=n\ln^\t n$, $\t>1$, or $k_n=n^2$.

Denote
\[
W_{n,k}=\sum_{j=-k_n+2}^{k_n} \Big(\D^{(2)} X_{s^n_j+t^n_k}\Big) =\sum_{j=-k_n+2}^{k_n}\left(X_{s^n_j+t^n_k}-2X_{s^n_{j-1}+t^n_k} +X_{s^n_{j-2}+t^n_k}\right)^2,
\]
where $1\le k\le n-1$ and $s^n_j=\frac{j}{m_n}T$.

\begin{theorem}\label{t:main} Assume that solution of Eq. (\ref{e:diflygt}) satisfies hypotheses $\mathbf{(H)}$ and $\mathbf{(H_1)}$. Moreover, let the
function $g$ is Lipschitz-continuous and there exists a random variable $\varsigma$ such that $\Prob(\varsigma<\infty)=1$ and inequality (\ref{i:boundness}) holds. Then
\begin{align*}
\widehat H^{(2)}_n&\tend{} H \quad\mbox{a.s.},\\
2\sqrt{n}\,\ln\frac{n}{T}\,(\widehat H^{(2)}_n-H)&\tendsd
N(0;\sigma^2_H)\qquad\mbox{for}\ H\in(1/2,1),
\end{align*}
where
\[
\widehat H_n^{(2)}=\frac{1}{2}+\frac{1}{2\ln k_n}\ln\left(\frac 2n
\sum_{k=2}^n\frac{\big(\D^{(2)} X_{t^n_k}
\big)^2}{W_{n,k-1}}\right),
\]
and $\sigma^2_H$ is a known variance defined in Subsection
\ref{s:fbm}.
\end{theorem}

\section{Preliminaries}\label{s:preliminaries}

\subsection{Several results on fBm}\label{s:fbm}

Recall that fBm $(B^H_t)_{t\ge0}$ with the Hurst index $H\in (0,1)$ is a
real-valued continuous centered Gaussian process with covariance
given by \[ \E(B^H_t B^H_s)
=\frac12\big(s^{2H}+t^{2H}-|t-s|^{2H}\big).
\]
For consideration of strong consistency and asymptotic normality of
the given estimators we need several facts regarding $B^H$.

\medskip\noindent\emph{\textbf{Limit results.}} Let
\[
V_{n,T}=\frac{n^{2H-1}}{T^{2H}(4-2^{2H})}\sum_{k=2}^n\big(\Delta^{(2)}
B^H_{t^n_k}\big)^2,\qquad H\neq \frac 12\,.
\]
Then (see \cite{Cohen-98}, \cite{IL}, \cite{begyn2})

\begin{align*}
    V_{n,T}&\tends{n\to\infty} 1\quad\text{ a.s.,} \\
   \sqrt{n}\left(V_{n,T}-1\right)&\tendsd\mathcal{N}\left(0,\sigma^2_H
   \right),
\end{align*}
where
\[
\sigma^2_H=2\bigg(1+2\sum_{j=1}^{\infty}\rho_H^2(j)\bigg),\quad
\rho_H(j)=-\frac{\abs{j-2}^{2H}-4\abs{j-1}^{2H}+6\abs{j}^{2H}-4\abs{j+1}^{2H} +\abs{j+2}^{2H} }{2(4-2^{2H})}\,.
\]

\medskip\noindent\emph{\textbf{H\"older-continuity of $B^H$.}} It is known that almost all sample paths of an fBm $B^H$ are locally H\"older of order strictly less than $H$, $H\in(0,1)$. To be more
precise, for all ƒ$0<\eps<H$ and $T>0$ there exists a nonnegative
random variable $G_{\eps,T}$ such that $\mathbb{E}(\vert
G_{\eps,T}\vert^p)<\infty$ for all $p\gs 1$, and
\begin{equation}\label{in:fBmpokyt}
\vert B^H_t -B^H_s\vert \ls G_{\eps,T}\vert t-s\vert^{H-\eps}\qquad
a.s.
\end{equation}
for all $s,t\in [0,T]$.

\subsection{Concentration inequality}

Let
\[
Y_{k,n}=\frac{n^H}{T^H\sqrt{4-2^{2H}} }\,\D^{(2)} B^H_{t^n_{k}},\qquad
d^{(2)n}_{jk}=\E Y_{j,n}Y_{k,n},\quad j,k=2,\dots,n.
\]
Note that $d^{(2)n}_{jk}=\rho_H(j-k)$. In the sequel we make use of the following
modified version of an inequality of concentration from \cite{bnp}.

\begin{lemma}\label{l:inequality} For all $z>0$ and any $H\in(0,1)$,
\[
\pr\bigg(\bigg\vert(n-1)^{-1/2}\sum_{k=2}^{n} (Y^2_{k,n} -1) \bigg\vert> z\bigg) \ls 2\,\exp\bigg(-\frac{z^2}{\frac{32}{3}(\frac{z}{\sqrt{n-1}}
+1)}\bigg).
\]
\end{lemma}
\proof Let $\kap=\sup_{H\in(0,1)}
\sum_{j\in\mathbb{Z}}\vert\rho_H(j)\vert$. Following an argument of
the paper \cite{bnp}, one gets bound
\begin{align*}
2(n-1)^{-1}\sum_{k,j=2}^{n} Y_{k,n} Y_{j,n}\, d^{(2)n}_{kj}
\ls&  2(n-1)^{-1}\sum_{k,j=2}^{n} \vert Y_{k,n}\vert\cdot \vert Y_{j,n}\vert\cdot \vert\rho_H(j-k)\vert\\
\ls& 2(n-1)^{-1}\sum_{k,j=2}^{n} Y_{k,n}^2\,
\vert\rho_H(j-k)\vert\\
\ls& 2(n-1)^{-1}\kap\sum_{k=2}^{n}Y_{k,n}^2=
2(n-1)^{-1}\kap\sum_{k=2}^{n} \big(Y^2_{k,n}-1\big)+ 2\kap\\
=&\frac{2\kap}{\sqrt{n-1}}\bigg(\frac{1}{\sqrt{n-1}} \sum_{k=2}^{n}
\big(Y^2_{k,n}-1\big)\bigg)+2\kap.
\end{align*}
Thus (see \cite{bnp}),
\[
\pr\bigg(\bigg\vert\frac{1}{\sqrt{n-1}}
\sum_{k=2}^{n} \big(Y^2_{k,n}-1\big)\bigg\vert> z\bigg) \ls
2\,\exp\bigg\{-\frac{z^2}{4\kap(\frac{z}{\sqrt{n-1}} +1)}\bigg\}.
\]
In paper \cite{bc} it was proved that
\[
\sum_{j\in\mathbb{Z}}\vert\rho_H(j)\vert=
\begin{cases}1+\frac{10-7\cdot 2^{2H}+2\cdot 3^{2H}}{4-2^{2H}}&\mbox{for}\quad H\ls 1/2,\\
1+\frac{4- 2^{2H}}{4-2^{2H}}=2,&\mbox{for}\quad H\gs 1/2,
\end{cases}
\]
and
\[
\kap=\sup_{H\in(0,1)}
\sum_{j\in\mathbb{Z}}\vert\rho_H(j)\vert=\lim_{H\to
0+}\sum_{j\in\mathbb{Z}}\vert\rho_H(j)\vert=\frac{8}{3} \,.
\]
This yields the required inequality.

\section{Proof of Theorem \ref{t:case}}\label{s:proofs1}

\noindent\emph{Proof of Theorem \ref{t:case}.} Observe first that
\[
\fy_{n,T}(x)=\Big(\frac{T}{n}\Big)^{2x}(4-2^{2x}),\qquad x\in(0,1),
\]
is continuous and strictly decreasing for $n>T$. Thus, it has an
inverse $\fy_{n,T}^{-1}$ for $n>T$. By hypothesis $\mathbf{(H_1)}$,
\begin{align*}
\frac{\fy_{n,T}(\widehat H^{(1)}_n)}{\fy_{n,T}(H)}=&\Big[\Big(\frac{T}{n}\Big)^{2H}(4-2^{2H})\Big]^{-1} \fy_{n,T}\bigg( \fy_{n,T}^{-1}\Bigg(\frac{1}{n}\sum_{i=2}^n\left(\frac{\Delta^{(2)}X_{t^n_i}}{g\big(X_{t^n_{i-1}}\big)}\right)^2\Bigg)\nonumber\\ =&\Big[\Big(\frac{T}{n}\Big)^{2H}(4-2^{2H})\Big]^{-1}\Bigg(\frac{1}{n}\sum_{i=2}^n \left(\frac{\Delta^{(2)}X_{t^n_i}}{g\big(X_{t^n_{i-1}}\big)}\right)^2\Bigg)\nonumber\\
=&\frac{n^{2H-1}}{T^{2H}(4-2^{2H})}\bigg(\sum_{i=2}^n \left(\Delta^{(2)}B^H_{{t^n_i}}\right)^2 +O_\omega\Big(\frac{1}{n^{3(H-\eps)-1}}\Big)\bigg)\nonumber\\
=&\frac{n^{2H-1}}{T^{2H}(4-2^{2H})}\sum_{i=2}^n
\left(\Delta^{(2)}B^H_{{t^n_i}}\right)^2+O_\omega\Big(\frac{1}{n^{H-3\eps}}\Big)=
V_{n,T}+O_\omega\Big(\frac{1}{n^{H-3\eps}}\Big)
\end{align*}
for $3\eps<H$. Therefore
\[
\frac{\fy_{n,T}(\widehat H^{(1)}_n)}{\fy_{n,T}(H)}\longrightarrow 1\qquad \mbox{a.s.\quad as}\ n\to\infty.
\]
Using the same arguments as in \cite{ksm} it is possible to prove  that the estimator $\widehat H^{(1)}_n$ is strongly consistent and asymptotically normal.

\section{Proof of the main Theorem}\label{s:main}

Before presenting the proof of this theorem, we  give two auxiliary
lemmas.

\begin{lemma}\label{l:trauka_prie_vid2}
Let
\[
    V_{n,T}(k) =\frac{m_n^{2H}}{2k_nT^{2H}(4-2^{2H})}\sum_{j=-k_n+2}^{k_n} \left(B^H_{s^n_j+t^n_k}-2B^H_{s^n_{j-1}+t^n_k} +B^H_{s^n_{j-2}+t^n_k}\right)^2,\qquad 1\ls k\ls n-1.
\]
The following relation holds:
\[
\max_{1\le k \le n-1}\Abs{V_{n,T}(k)-1}
=O_\omega\bigg(\sqrt{\frac{\ln n}{k_n}}\bigg).
\]
\end{lemma}
\noindent\emph{Proof.} By self similarity and stationarity of
increments of fBm,
\begin{align*}
    V_{n,T}(k)\stackrel{d}{=}&\frac{2^{2H-1}m_n^{2H}}{2^{2H}k_n(4-2^{2H})}\sum_{j=-k_n+2}^{k_n}\left(B^H_{\frac{j+k_n}{m_n}} -2B^H_{\frac{j+k_n-1}{m_n}}+B^H_{\frac{j+k_n-2}{m_n}}\right)^2\\
        \stackrel{d}{=}&\frac{(2m_n)^{2H-1}n}{(2n)^{2H}(4-2^{2H})}\sum_{j=-k_n+2}^{k_n}\left(B^H_{\frac{j}{k_n}+1}-2B^H_{\frac{j-1}{k_n}+1} +B^H_{\frac{j-2}{k_n}+1}\right)^2\\
   \stackrel{d}{=}&\frac{(2k_n)^{2H-1}}{2^{2H}(4-2^{2H})}\sum_{j=2}^{2k_n}\left(B^H_{\frac{j}{k_n}}-2B^H_{\frac{j-1}{k_n}} +B^H_{\frac{j-2}{k_n}}\right)^2\\
    \stackrel{d}{=}& \frac{(2k_n)^{2H-1}}{4-2^{2H}}\sum_{j=2}^{2k_n}\left(B^H_{\frac{j}{2k_n}}-2B^H_{\frac{j-1}{2k_n}}+B^H_{\frac{j-2}{2k_n}}\right)^2 \stackrel{d}{=}V_{2k_n,1}.
\end{align*}

Therefore,
\[
\Prob\left(\max_{1\le k \le n-1}\Abs{V_{n,T}(k)-1}>\de\right)\le
\sum_{k=1}^{ n-1} \Prob\left(\Abs{V_{n,T}(k)-1}>\de\right) \le
n\Prob\left(\abs{V_{2k_n,1}-1}>\de\right)\qquad \mbox{\rm for all}\
\de>0.
\]
Put $\widehat{V}_{n,T}=\frac{n}{n-1}V_{n,T}$. Note that
\[
\big\vert V_{2k_n,1}-1\big\vert\le \big\vert\widehat{V}_{2k_n,1}-1\big\vert+\frac{1}{2k_n}\,.
\]

Let $(\de_n)$ be a sequence of positive numbers such that
$\de_n\downarrow 0$ as $n\to\infty$ and $k_n^{-1}<\de_n$. By Lemma
\ref{l:inequality},
\begin{align*}
    \Prob\left(\abs{V_{2k_n,1}-1}>2\de_n\right)\le&
    \Prob\left(\abs{\widehat{V}_{2k_n,1}-1}+\frac{1}{2k_n}>2\de_n\right) \le \Prob\left(\abs{\widehat{V}_{2k_n,1}-1}>\de_n\right)\\ =&\Prob\left(\sqrt{2k_n-1}\,\abs{\widehat{V}_{2k_n,1}-1}>\de_n\sqrt{2k_n-1}\right)
    \le 2\,\exp\bigg\{-\frac{\de^2_n(2k_n-1)}{\frac{32}{3}(\de_n +1)}\bigg\}  \,.
\end{align*}
Set  $\de_n=\sqrt{a\frac{\ln n}{2k_n-1}}$. Since $k_n\gs n\ln n$, then
\[
\Prob\left(\abs{V_{2k_n,1}-1}>2\de_n\right)\le 2\,\exp\left\{-\frac{3a\ln n}{32\left(\sqrt{a\frac{\ln n}{2n\ln n-1}}+1\right)}\right\}.
\]
If $a\gs 3$ and $n\gs 2$, then $\de_n>k_n^{-1}$. Moreover,
$\Prob\left(\abs{V_{2k_n,1}-1}>2\de_n\right)\le 2 n^{-3}$ for $a>
32$ and $n$ large enough. Therefore series
$\sum\limits_n\Prob(\max_{1\le k \le n-1}\Abs{V_{n,T}(k)-1}>\de_n)$
converges and by the Borel-Cantelli lema $\max_{1\le k \le
n-1}\Abs{V_{n,T}(k)-1}\tends{n\to\infty} 0$ a.s.

\begin{lemma}\label{l:nario_asimp}
Assume that function $g$ is Lipschitz-continuous. If
$\eps<(H-1/2)/3$, then for each $k=1,\dots,n-1$
\[
W_{n,k}=g^2\left(X_{t^n_k}\right)\frac{T^{2H}(4-2^{2H})}{(2k_n)^{2H-1}}
+ O_\omega\left(\frac{\sqrt{k_n\ln n}}{m_n^{2H}}\right)
+O_\omega\left(\frac{k_n}{n^{H-\eps}m_n^{2(H-\eps)}}\right).
\]
\end{lemma}
\noindent\emph{Proof. Step 1.} By hypothesis $\mathbf{(H_1)}$,
\begin{equation*}
 \Delta^{(2)}
    X_{s^n_j+t^n_k}=g\big(X_{s^n_{j-1}+t^n_k}\big)\Delta^{(2)}
    B^{H}_{s^n_j+t^n_k}+O_\omega\left(\frac{1}{m_n^{2(H-\eps)}}\right),\qquad j=-k_n+2,\dots, k_n.
\end{equation*}
Next, note that
\begin{align*}
\abs{g^2(X_t)-g^2(X_s)}=&\abs{g(X_t)-g(X_s)}\abs{g(X_t)+g(X_s)}\le
L\abs{X_t-X_s}\abs{g(X_s)+g(X_t)}\\
\ls& 2L\sup_{u\in[0,T]}\abs{g(X_u)}
\abs{X_t-X_s}
\end{align*}
where $t>s$ and $L$ is Lipschitz constant. Thus, hypothesis
$\mathbf{(H)}$ with a.s. continuity of $t\mapsto g(X_t)$ lead to
\[
g^2\left(X_{s^n_{j-1}+t^n_k}\right) -g^2\left(X_{t^n_k}\right)=O_\omega\bigg(\frac{1}{n^{H-\eps}}\bigg).
\]

\noindent\emph{Step 2.} Assume that $\eps<(H-1/2)/3$. \emph{Step 1}
and Lemma \ref{l:trauka_prie_vid2} yield
\begin{align*}
    W_{n,k}=&\sum_{j=-k_n+2}^{k_n}\left(\Delta^{(2)}X_{s^n_j+t^n_k}\right)^2=
    \sum_{j=-k_n+2}^{k_n} g^2\left(X_{s^n_{j-1}+t^n_k}\right)\left(\Delta^{(2)}B^H_{s^n_{j}+t^n_k}\right)^2 +O_\omega\left(\frac{k_n}{m_n^{3(H-\eps)}}\right)\\
    =&g^2\left(X_{t^n_k}\right)\sum_{j=-k_n+2}^{k_n}\left(\Delta^{(2)}B^{H}_{s^n_{j}+t^n_k}\right)^2+
    \sum_{j=-k_n+2}^{k_n}\left(g^2\left(X_{s^n_{j-1}+t^n_k}\right) -g^2\left(X_{t^n_k}\right)\right) \left(\Delta^{(2)}B^{H}_{s^n_{j}+t^n_k}\right)^2\\ &+O_\omega\left(\frac{k_n}{m_n^{3(H-\eps)}}\right)\\
    =& g^2\left(X_{t^n_k}\right)\frac{2k_nT^{2H}(4-2^{2H})}{m_n^{2H}} \,V_{n,T}(k)+O_\omega\left(\frac{k_n}{n^{H-\eps}m_n^{2(H-\eps)}}\right) +O_\omega\left(\frac{k_n}{m_n^{3(H-\eps)}}\right)\\
    =& g^2\left(X_{t^n_k}\right)\frac{2k_nT^{2H}(4-2^{2H})}{m_n^{2H}} + O_\omega\left(\frac{\sqrt{k_n\ln n}}{m_n^{2H}}\right)
+O_\omega\left(\frac{k_n}{n^{H-\eps}m_n^{2(H-\eps)}}\right).
\end{align*}
Consequently the proof of lemma is completed.

\smallskip\emph{Proof of Theorem \ref{t:main}}. Put
\[
S_{n,T}:=\frac{2}{nk_n^{2H-1}} \sum_{k=2}^n\frac{\big(\D^{(2)} X_{t^n_k}\big)^2}{W_{n,k-1}}\,.
\]
It follows from (\ref{e:nario_asimp_H1})--(\ref{i:boundness}) and
Lemma \ref{l:nario_asimp} that
\begin{align*}
S_{n,T}=&\frac{2}{nk_n^{2H-1}}\sum_{k=2}^n\frac{g^2\left(X_{t^n_k}\right) \big(\Delta^{(2)}
B^{H}_{t^n_k}\big)^2+O_\omega(n^{-3(H-\eps)})}{g^2 \left(X_{t^n_k}\right)\frac{2k_nT^{2H}(4-2^{2H})}{m_n^{2H}} +O_\omega\left(\frac{\sqrt{k_n\ln n}}{m_n^{2H}}\right)+O_\omega\left(\frac{k_n}{n^{H-\eps}m_n^{2(H-\eps)}}\right)}\\
=&\frac{2}{nk_n^{2H-1}}\sum_{k=2}^n\frac{\big(\Delta^{(2)}
B^{H}_{t^n_k}\big)^2+O_\omega\left(n^{-3(H-\eps)}\right)}{\frac{2k_nT^{2H}(4-2^{2H})}{m_n^{2H}} +O_\omega\left(\frac{\sqrt{k_n\ln n}}{m_n^{2H}}\,\right)+O_\omega\left(\frac{k_n}{n^{H-\eps}m_n^{2(H-\eps)}}\right)}\\
=&\frac{m_n^{2H}}{nk_n^{2H}T^{2H}(4-2^{2H})} \sum_{k=2}^n\frac{\big(\Delta^{(2)}
B^{H}_{t^n_k}\big)^2+O_\omega\left(n^{-3(H-\eps)}\right)}{1 +O_\omega\left(\sqrt{\frac{\ln n}{k_n}}\,\right)+O_\omega\left(\frac{m_n^{2\eps}}{n^{H-\eps}}\right)}\\
=&\frac{V_{n,T}+O_\omega\big(n^{-(H-3\eps)}\big)}{1 +O_\omega\left(\sqrt{\frac{\ln n}{k_n}}\,\right)+O_\omega\left(\frac{m_n^{2\eps}}{n^{H-\eps}}\right)}\,.
\end{align*}

Term $O_\omega\left(\frac{m_n^{2\eps}}{n^{H-\eps}}\right)$ vanishes
provided $\eps>0$ is small enough. Convergence $V_{n,T}\tendsb{}1$
implies that $S_{n,T}\tendsb{}1$. Hence
\[
H^{(2)}_n=H+\frac{\ln S_{n,T}}{2\ln
k_n}\underset{n\to\infty}{\tendsb}H.
\]

To prove asymptotic normality of the estimator $H^{(2)}_n$ observe
that
\begin{align*}
\sqrt{n}\,\big(S_{n,T}-1\big) =&\sqrt{n}\,\Bigg(\frac{V_{n,T}-1+O_\omega\left(\sqrt{\frac{\ln n}{k_n}}\right)+O_\omega\left(\frac{m_n^{2\eps}}{n^{H-\eps}}\right)  +O_\omega(n^{-(H-3\eps)})}{1+O_\omega\left(\sqrt{\frac{\ln n}{k_n}}\right)+O_\omega\left(\frac{m_n^{2\eps}}{n^{H-\eps}}\right)} \Bigg)\\
=&\frac{\sqrt{n}(V_{n,T}-1)}{1+O_\omega\left(\sqrt{\frac{\ln
n}{k_n}}\right)
+O_\omega\left(\frac{m_n^{2\eps}}{n^{H-\eps}}\right)}\\
&+\frac{O_\omega\left(\sqrt{\frac{n\ln n}{k_n}}\,\right)
+O_\omega\left(\frac{m_n^{2\eps}}{n^{H-1/2-\eps}}\right)
+O_\omega\left(n^{-(H-1/2-3\eps)}\right)
}{1+O_\omega\left(\sqrt{\frac{\ln n}{k_n}}\right)
+O_\omega\left(\frac{m_n^{2\eps}}{n^{H-\eps}}\right)}\tendsd{}N(0,\s_H^2)
\end{align*}
for $\eps>0$ small enough. Now apply Slutsky's theorem and limit
results of Section \ref{s:fbm}.

\section{Examples}\label{s:examples}
As mentioned previously, in this section we present two examples of
applications of the obtained results. The first one deals with a
general form of equation \eqref{e:diflygt} and relies on certain
restrictions on functions $f$ and $g$. The second one describes a
particular model which formally does not fit into the scope of the
first one.

\subsection{Example 1}\label{ss:example1}

In order to present an example, we need several facts on variation.
To make the paper more self-contained and the structure clearer, the
mentioned facts are briefly reminded in subsection
\ref{ss:variation}. For details we refer the reader to
\cite{DudleyNorvaisa-10}.

\subsubsection{Variation}\label{ss:variation}

Fix $p>0$ and $-\infty<a<b<\infty$. Let
$\varkappa=\{\{x_0,\dots,x_n\}\mid a=x_0<\dots<x_n=b,n\ge1\}$
denotes a set of all possible partitions of $[a,b]$. For any
$f:[a,b]\to\mathbb{R}$ define
\begin{gather*}
    v_{p}(f;[a,b])=\sup_{\varkappa}\sum_{k=1}^{n}\abs{f(x_k)-f(x_{k-1})}^{p},\qquad
    V_{p}(f;[a,b])=v_{p}^{1/p}(f;[a,b]),\\
    \mathcal{W}_p([a,b])=\{f:[a,b]\to\mathbb{R}\mid
    v_{p}(f;[a,b])<\infty\},\quad
    C\mathcal{W}_p([a,b])=\{f\in\mathcal{W}_p([a,b])\mid f\text{ is
    continuous}\}.
\end{gather*}
Recall that $v_p$ is called $p$-variation of $f$ on $[a,b]$ and any
$f$ in $\mathcal{W}_p([a,b])$ is said to have bounded $p$-variation
on $[a,b]$. For short we omit an interval $[a,b]$ in the notations
introduced above whenever there is no ambiguity. Below we list
several facts used further on.

\begin{itemize}
  \item $f\mapsto V_p(f)$ is a seminorm on $\mathcal{W}_p$; $V_p(f)=0$ if and only if $f$ is a constant.
  \item $f\in \mathcal{W}_p\imply
  \sup_{x\in[a,b]}\abs{f(x)}<\infty$.
  \item $f,g\in \mathcal{W}_p\imply fg\in \mathcal{W}_p$.
  \item $q>p\gs 1\imply \mathcal{W}_{p}\subset \mathcal{W}_{q}$.
  \item Let $f\in \mathcal{W}_q$, $h\in \mathcal{W}_p$ with
        $p,q\in(0,\infty)$ such that  $1/p+\allowbreak 1/q>1.$ Then an integral
        $\int_a^b f\,\mathrm{d}h$ exists as the Riemann--Stieltjes  integral
        provided $f$ and $h$ have no common discontinuities. If the integral
        exists, the Love--Young inequality
        \[
        \Bigg\vert \int\limits_a^bf\,\dr h-f(y)\big[ h(b)-h(a)
        \big]\Bigg\vert \ls C_{p,q} V_q\big(f\big)V_p\big(h\big)
        \]
        holds for all $y\in [a,b]$, where $C_{p,q}=\zeta(p^{-1}+q^{-1})$ and
        $\zeta(s)=\sum_{n\gs 1} n^{-s}$. Moreover,
       $$
       V_p\Bigg (\int\limits_a^{\bt} f \,\dr h; [a, b] \Bigg) \leq C_{p, q}
       V_{q,\infty}\big(f\big) V_p\big(h\big),
       $$
       where $V_{q,\infty}(f)=V_q(f)+\sup_{x\in[a, b]}\vert
       f(x)\vert$. Also note that $V_{q,\infty}$ is a norm on ${\cal
       W}_q$, $q\gs 1$.
\end{itemize}

\begin{remark}
The left-hand side of \eqref{in:fBmpokyt} can be replaced by
$V_{H_\eps}(B^{H};[s,t])$,   $\forall\eps\in(0,H-\frac{1}{2})$, i.e.
\begin{equation}\label{e:nel_fBm_variacijai}
    V_{H_\eps}(B^{H};[s,t])\le
    G_{\eps,T}\abs{t-s}^{H-\eps}\text{ a.s.}
\end{equation}
with $G_{\eps,T}$ of \eqref{in:fBmpokyt}.
\end{remark}

\subsubsection{Assumptions and properties of solution of SDE}\label{ss:solution}

Let $(B^H_t)_{t\in[0,T]}$, $H\in(1/2,1)$, be a fixed fBm defined on some probability space $(\Omega,\,{\mathcal F},\,{\bf P},\,\mathbb{F})$. Let $\alpha\in\left(\frac{1}{H}-1;1\right]$,
$\mathcal{C}^{1+\a}(\mathbb{R})=\{h:\Rd\to\Rd\mid h^\prime\text{ exists
and }\sup_x\abs{h^\prime(x)}+\sup_{x\neq
y}\frac{\abs{h^\prime(x)-h^\prime(y)}}{\abs{x-y}^\alpha}<\infty\}$. Assume that $f$ is Lipschitz and
$g\in\mathcal{C}^{1+\a}(\mathbb{R})$. In such case there exists a
unique solution of \eqref{e:diflygt} having the following
properties: i) $(X_t)_{t\in[0,T]}$ is $\mathbb{F}$ adapted and
almost all sample paths are continuous; ii) $X_0=\xi$ a.s.; iii)
$\forall p>\frac{1}{H},\,\Prob(V_p(X;[0,T])<\infty)=1$ (see
\cite{du}, \cite{ly}, \cite{ly1} and  \cite{kk}).

\begin{lemma} Let $X$ satisfies
(\ref{e:diflygt}), $\eps\in(0,H-\frac{1}{2})$.
There exists a.s. finite r.v. $L_{\eps,T}$  such that
\begin{equation}\label{e:nel_variacijai}
V_{H_\eps}\big(X;[s,t]\big)\ls L_{\eps,T}\,(t-s)^{(H-\eps)}, \qquad
0\le s<t\le T.
\end{equation}

\end{lemma}
\noindent\emph{Proof.} Since $V_p$
is seminorm and non-increasing function of $p\gs 1$, inequalities of
subsection \ref{ss:variation} give bound
\begin{align*}
V_{H_\eps}\big(X;[s, t]\big)\ls& \int_s^t\abs{f(X_u)}\,du +V_p\bigg(\int_s^{\bt} g(X_u)\,d B^{H}_u;[s, t]\bigg)\\
\ls&\sup_{u\in [0,T]}\abs{f(X_u)}(t-s) +
C_{H_\eps,H_\eps}\left(V_{H_\eps}\big(g\circ X_{\cdot};[s,
t]\big)+\sup_{u\in
[0,T]}\abs{g(X_u)}\right)V_{H_\eps}\big(B^H;[s,t]\big)\\
\ls& L_{\eps,T}(t-s)^{H-\eps}
\end{align*}
with $L_{\eps,T}=\sup_{u\in [0,T]}\abs{f(X_u)}T^{1-H+\eps} +
C_{H_\eps,H_\eps}\left(V_{H_\eps}\big(g\circ X_{\cdot};[s,
t]\big)+\sup_{u\in [0,T]}\abs{g(X_u)}\right)G_{\eps,T}$ and
$G_{\eps,T}$ of \eqref{e:nel_fBm_variacijai}. A. s. continuity of
$t\mapsto X_t$ together with continuity of $f$ and $g$ implie a.s.
boundedness of $\sup_{u\in [0,T]}\abs{f(X_u)}$ and $\sup_{u\in
[0,T]}\abs{g(X_u)}$. It is not difficult to show that
$V_{H_\eps}\big(g\circ X_{\cdot};[s, t]\big)\le \sup_{u}\abs{g^\prime(u)}\,V_{H_\eps}\big(X;[s, t]\big)$. Hence, continuity
of $g^\prime$ and a.s. boundedness of $V_{H_\eps}\big(X;[s, t]\big)$
guarantees a.s. boundedness of $V_{H_\eps}\big(g\circ X_{\cdot};[s,
t]\big)$ along with that of $L_{\eps,T}$.

\begin{lemma}\label{l:nario_asimp0} Let $X$ satisfies (\ref{e:diflygt}), $\eps\in(0,H-\frac{1}{2})$,
and conditions stated above are true. Then the following relations
hold:
\begin{align}
\Delta X_{\tau^n_k}=& O_\omega\big(d_n^{H-\eps}\big),\qquad k=1,\dots, i_n,\label{e:nario_asimp_0}\\
 \Delta^{(2)} X_{\tau^n_k} =&g(X_{\tau^n_{k-1}})\Delta^{(2)}
    B^{H}_{\tau_k}+O_\omega\big(d_n^{2(H-\eps)}\big),\qquad k=2,\dots,
    i_n,\label{e:nario_asimp_1}
\end{align}
where $d_n=\max_{1\ls k\ls i_n}(\tau^n_k-\tau^n_{k-1})$.
\end{lemma}
\noindent\emph{Proof.} Let a sample path $t\mapsto X_t$ be
continuous. We first prove (\ref{e:nario_asimp_0}). Note that
\begin{equation*}
\Delta X_{\tau^n_k}=X_{\tau^n_{k}}-X_{\tau^n_{k-1}}=
\int_{{\tau^n_{k-1}}}^{\tau^n_{k}} f(X_s) ds
+\int_{{\tau^n_{k-1}}}^{\tau^n_{k}}
\big[g(X_s)-g(X_{\tau^n_{k-1}})\big]dB^H_s +
g(X_{\tau^n_{k-1}})\Delta B^H_{\tau^n_k}.
\end{equation*}
An application of inequalities
\eqref{e:nel_fBm_variacijai}--\eqref{e:nel_variacijai} together with
continuity of $f$, $g,g^\prime$ and mean value theorem yield
\begin{align}
& \int_{\tau^n_{k-1}}^{\tau^n_{k}} \left| f(X_s)\right|ds\ls \sup_{u\in[0,T]}\abs{f(X_u)} \, (\tau^n_k-\tau^n_{k-1})=O_\omega(d_n), \nonumber \\
&\left|\int_{\tau^n_{k-1}}^{\tau^n_{k}}
[g(X_s)-g(X_{\tau^n_{k-1}})]\,dB^H_s\right| \ls
\norm{g^\prime}_\infty
C_{{H_\eps},{H_\eps}}V_{H_\eps}\big(X;[\tau^n_{k-1},\tau^n_{k}]\big)
V_{H_\eps}\big(B^H;[\tau^n_{k-1},\tau^n_{k}]\big)\nonumber\\
&\qquad=O_\omega\big(d_n^{2(H-\eps)}\big),\label{e:nel_integralui2}\\
&\abs{g(X_{\tau^n_{k}})\Delta
B^H_{\tau^n_{k}}}\le\sup_{u\in[0,T]}\abs{g(X_u)}
G_{\eps,T}(\tau^n_k-\tau^n_{k-1})^{H-\eps}=O_{\omega} (d_n^{H-\eps}).\nonumber
\end{align}
Therefore \eqref{e:nario_asimp_0} holds.

Next we prove (\ref{e:nario_asimp_1}). Since
\[
\Delta^{(2)}X_{\tau^n_k}=\int_{\tau^n_{k-1}}^{\tau^n_k}f(X_s)\d
    s -\int_{\tau^n_{k-2}}^{\tau^n_{k-1}}f(X_s)\d s + \int_{\tau^n_{k-1}}^{\tau^n_k}g(X_s)\d B^{H}_
    s-\int_{\tau^n_{k-2}}^{\tau^n_{k-1}}g(X_s)\d B^{H}_s,
\]
and
\begin{align*}
    &\int_{\tau^n_{k-1}}^{\tau^n_k}g(X_s)\d
    B^{H}_s-\int_{\tau^n_{k-2}}^{\tau^n_{k-1}}g(X_s)\d B^{H}_s
    =\int_{\tau^n_{k-1}}^{\tau^n_k}[g(X_s)-g(X_{\tau^n_{k-1}})]\d B^{H}_s\\
    &\quad-\int_{\tau^n_{k-2}}^{\tau^n_{k-1}}[g(X_s)-g(X_{\tau^n_{k-1}})]\d B^{H}_s
    +g(X_{\tau^n_{k-1}})\left(\Delta B^H_{\tau_k^n}-\Delta B^H_{\tau^n_{k-1}}\right)\\
    &\quad\stackrel{\eqref{e:nel_integralui2}}{=}
    O_{\omega}(d_n^{2(H-\eps)})+g(X_{\tau^n_{k-1}})\Delta^{(2)} B^H_{\tau_k^n},\qquad
    k=2,\dots,i_n,
\end{align*}
then by mean value theorem with $K$ equal to Lipschitz constant of $f$,
\begin{align*}\label{e:nel_integralui3}
    &\left|\int_{\tau^n_{k-1}}^{\tau^n_k}f(X_s)\d s-\int_{\tau^n_{k-2}}^{\tau^n_{k-1}}f(X_s)\d
    s\right|\le \int_{\tau^n_{k-1}}^{\tau^n_k}\left|f(X_s)-f(X_{\tau^n_{k-1}})\right|\d s+\int_{\tau^n_{k-2}}^{\tau^n_{k-1}}\left|f(X_{\tau^n_{k-1}})-f(X_s)\right|\d s\nonumber\\
    &\quad\le 2K d_n\max_{1\le k\le i_n}\sup_{\tau^n_{k-1}\le s\le
    \tau^n_k}\abs{X_s-X_{\tau^n_{k-1}}}
    \stackrel{\eqref{e:nario_asimp_0}}{=}O_{\omega}\big(d_n^{1+H-\eps}\big).
\end{align*}
We conclude that \eqref{e:nario_asimp_1} also holds.

\subsubsection{Application}

Results of subsection \ref{ss:solution} imply proposition
constituting the basis of the first example.

\begin{proposition}
Assume that a model defined by \eqref{e:diflygt} satisfies
conditions of subsection \ref{ss:solution}. If in addition
\eqref{i:boundness} holds, then Theorems \ref{t:case}--\ref{t:main}
apply.
\end{proposition}

\subsection{Example 2}

Consider the Verhulst equation
\[
X_t=\xi+\int_0^t (\l X_s-X_s^2)\, dt + \s \int_0^t X_s\, d B^H_t, \quad \xi>0,\quad  t\in[0,T].
\]
$f(x)=\l x-x^2$ is not Lipschitz. So we can't use results of section
\ref{ss:example1}.  It was proved in \cite{ksm} that this equation
admits an explicit solution
\[
X_t= \frac{\xi \exp\{\l t + \s B^H_t\}}{1+ \xi\int_0^t \exp\{\l s +\s B^H_s\}ds}, \quad  t\in[0,T].
\]

\noindent Moreover, from \cite{ksm} one can get
\begin{align*}
\Delta X_{t^n_k}=&X_{t^n_{k-1}}\bigg(\l \,\frac{T}{n} + \s \Delta B^H_{t^n_k} +
O_\omega\bigg(\frac{1}{n}\bigg)\bigg),\qquad k=1,\dots, n,\\
\Delta^{(2)} X_{t^n_k} =& X_{t^n_{k-1}}\left(\s\Delta^{(2)} B^{H}_{t^n_k}+ O_\omega\left(\frac{1}{n^{2(H-\eps)}}\right)\right),\qquad k=2,\dots, n.
\end{align*}
Finally, explicit form of $X_t$ implies existence of an a.s. finite
r.v. $\varsigma$ such that $ \sup_{t\in[0,T]}\frac{1}{\vert
X_t\vert}\le \varsigma\quad\mbox{a.s.}$ Thus, one can use Theorem
\ref{t:case} if constant $\s$ is known and Theorem \ref{t:main} in
general situation.


\end{document}